\documentclass[journal,10pt,onecolumn]{IEEEtran}
\usepackage{amsmath}
\usepackage{amsfonts}
\usepackage{amsthm}
\usepackage{dsfont}
\usepackage{epsfig}
\usepackage{upgreek}
\usepackage{graphicx}

\usepackage{setspace}

\theoremstyle{definition}

\newcommand{\Jm}{\mbox{${\bf J}$}}

\newcommand{\Imat}{\mbox{${\bf I}$}}

\newcommand{\Amn}{\mbox{${\bf A}$}}
\newcommand{\Amni}{\mbox{${\bf A}^{\!\mbox{\tiny -1}}$}}

\newcommand{\Sm}{\mbox{${\bf S}$}}
\newcommand{\Smi}{\mbox{${\bf S}^{\!\mbox{\tiny -1}}$}}

\newcommand{\Qm}{\mbox{${\bf Q}$}}
\newcommand{\Qmt}{\mbox{${\bf Q}^{\mbox{\tiny{T}}}$}}

\newcommand{\Pm}{\mbox{${\bf P}$}}

\newcommand{\Vm}{\mbox{${\bf V}$}}
\newcommand{\Vmi}{\mbox{${\bf V}^*$}}
\newcommand{\Xm}{\mbox{${\bf X}$}}

\newcommand{\Ym}{\mbox{${\bf Y}$}}

\newcommand{\Mm}{\mbox{${\bf M}$}}
\newcommand{\Dm}{\mbox{${\bf D}$}}
\newcommand{\Dmi}{\mbox{${\bf D}^{\mbox{\tiny -1}}$}}

\newcommand{\Em}{\mbox{${\bf E}$}}
\newcommand{\Emi}{\mbox{${\bf E}^{\mbox{\tiny -1}}$}}

\newcommand{\dinv}[1]{{#1}^{\mbox{\tiny -D}}}

\newcommand{\ginv}[1]{{#1}^{\mbox{\tiny -U}}}
\newcommand{\pinv}[1]{{#1}^{\mbox{\tiny -P}}}
\newcommand{\sinv}[1]{{#1}^{\mbox{\tiny -S}}}
\newcommand{\inv}[1]{{#1}^{\mbox{\tiny -1}}}

\newcommand{\Um}{\mbox{${\bf U}$}}
\newcommand{\Umi}{\mbox{${\bf U}^*$}}

\newcommand{\Atinv}{\mbox{${\bf A}^{\overset{\sim}{\!\mbox{\tiny -1}}}$}}
\newcommand{\Adinv}{\mbox{${\bf A}^{\!\mbox{\tiny -D}}$}}
\newcommand{\Apinv}{\mbox{${\bf A}^{\!\mbox{\tiny -P}}$}}
\newcommand{\Aginv}{\mbox{${\bf A}^{\!\mbox{\tiny -U}}$}}
\newcommand{\Asinv}{\mbox{${\bf A}^{\!\mbox{\tiny -S}}$}}

\newcommand{\tinv}[1]{\mbox{${#1}^{\overset{\sim}{\mbox{\tiny -1}}}$}}

\newcommand{\yv}{\mbox{$\hat{\bf y}$}}
\newcommand{\av}{\mbox{$\hat{\bf \uptheta}$}}
\newcommand{\yyv}{\mbox{$\hat{\bf y}'$}}
\newcommand{\aav}{\mbox{$\hat{\bf \uptheta}'$}}

\begin{document}

\title{{A Rank-Preserving Generalized Matrix Inverse for\\
           Consistency with Respect to Similarity} }       
\author{
\IEEEauthorblockN{{Jeffrey Uhlmann}}\\
\IEEEauthorblockA{\small University of Missouri-Columbia\\
201 EBW, Columbia, MO 65211\\
Email: uhlmannj@missouri.edu}}
\date{}          
\maketitle
\thispagestyle{empty}


\begin{abstract}
There has recently been renewed recognition of the need to understand the consistency properties
that must be preserved when a generalized matrix inverse is required. The most widely known generalized inverse,
the Moore-Penrose pseudoinverse, provides consistency with respect to orthonormal transformations
(e.g., rotations of a coordinate frame), and a recently derived inverse provides consistency with 
respect to diagonal transformations (e.g., a change of units on state variables). Another well-known
and theoretically important generalized inverse is the Drazin inverse, which preserves consistency
with respect to similarity transformations. In this paper we note a limitation of the Drazin inverse
is that it does not generally preserve the rank of the linear system of interest. We then 
introduce an alternative generalized inverse that both preserves rank and 
provides consistency with respect to similarity transformations. We discuss practical
implementation considerations and demonstrate with an example.\\ 

\begin{footnotesize}
\noindent {\bf Keywords}: {\sf\scriptsize Control Systems, Drazin Inverse, Generalized Matrix Inverse, 
Inverse Problems, Linear Estimation, Linear Systems, Matrix Analysis,  Moore-Penrose Pseudoinverse, 
Relative Gain Array, RGA, Similarity Transformations, Singular Value Decomposition, SVD,
Stability of Linear Systems, System Design, UC inverse, Unit Consistency.}  
\end{footnotesize}
\end{abstract}

\section{Introduction}

Many analytical and practical problems require a solution to an underdetermined
or overdetermined system of equations in the form of a generalized matrix 
inverse of a singular matrix $\Amn$ as $\Atinv$ satisfying
\begin{eqnarray}
      \Amn\Atinv\Amn & = & \Amn,\\
      \Atinv\Amn\Atinv\!\! & = & \Atinv.
\end{eqnarray}
In addition to these properties, a practical application may require
the system to be consistent with respect to a general class of 
transformations. For example, if the behavior of the system should
be invariant with respect to arbitrary rotations of the Cartesian
coordinate frame in which it is defined then the appropriate 
generalized inverse should also satisfy
\begin{equation}
     \tinv{(\Um\Amn\Vm)} ~ = ~ \Vmi\Atinv\Umi \label{udist}
\end{equation}
where $\Um$ and $\Vm$ are arbitrary rotation matrices (or, more generally,
arbitrary orthonormal/unitary matrices) such that $\Um\Umi=\Imat$, where
$\Umi$ is the conjugate-transpose of $\Um$. This condition is satisfied
by the Moore-Penrose inverse (MP inverse), which we will denote as
$\Apinv$. In other applications the behavior of the system of interest is 
expected to be invariant with respect to changes of units defined for 
the its key state variables/parameters. For example, if all variables
are defined consistently in imperial units then a change to metric
units should have no effect on the behavior of the system, i.e., the
system should be invariant to the arbitrary choice of units. This condition can
be expressed as 
\begin{eqnarray}
      \tinv{(\Dm\Amn\Em)}  & = & \Emi\Atinv\Dmi
\end{eqnarray}
where $\Dm$ and $\Em$ are arbitrary nonsingular diagonal
matrices. This condition is satisfied by a recently-derived
unit-consistent (UC) generalized inverse~\cite{uhlmann0}.

An unfortunate fact is that in many scientific and engineering
domains the MP inverse is widely used 
as the default generalized inverse without consideration
for the properties it preserves -- {\em and does not preserve}.
Examples include some forms of Channel Estimation and 
some applications of the Relative Gain Array (RGA).
In particular, the RGA is used to
measure process variable interactions defined by a matrix $\Amn$
and is defined as~\cite{bristol}
\begin{equation}
   \mbox{RGA}(\Amn) ~ \doteq ~ \Amn \circ (\inv{\Amn})^{\mbox{\scriptsize T}}.
\end{equation}
where $\circ$ represents the elementwise Hadamard matrix product.
When $\Amn$ is singular (or non-square) the RGA is often 
evaluated in practice as~\cite{changyu}
\begin{equation}
    \Amn \circ (\Apinv)^{\mbox{\scriptsize T}}
\end{equation}
without regard for whether the relative magnitudes of interactions
should be invariant to the choice of units on input/output
variables. In many applications involving singular system
matrices (likely most) unit consistency is needed 
and thus the RGA should be evaluated using the UC inverse as
\begin{equation}
    \Amn \circ (\Aginv)^{\mbox{\scriptsize T}}
\end{equation}
Failure to do so can lead to erroneous pairings of controlled 
and manipulated variables due to inappropriate dependencies
on the choice of units. 

In other applications there may be need for an
$n\times n$ generalized
inverse to be consistent with respect to an arbitrary matrix
similarity transformation:
\begin{equation}
     \tinv{(\Sm\Amn\Smi)} ~ = ~ \Sm\Atinv\Smi \label{dinvcnst}
\end{equation}
The need for this condition can be seen from a standard linear model
\begin{equation}
     \yv ~ = ~ \Amn\cdot\av \label{linmod}
\end{equation}
where the objective is to identify a vector $\av$ of parameter values
satisfying the above equation for a data matrix $\Amn$ and
a known/desired state vector $\yv$ defined in the same state
space as $\av$. If $\Amn$ is nonsingular
then there exists a unique $\Amni$ which gives the solution
\begin{equation}
     \av ~ = ~ \Amni\cdot\yv
\end{equation}
If $\Amn$ is singular, however, then the
Moore-Penrose inverse could be applied as
\begin{equation}
     \av ~ = ~ \Apinv\cdot\yv
\end{equation}
or the UC inverse could be applied as
\begin{equation}
     \av ~ = ~ \Aginv\cdot\yv
\end{equation}
to obtain a solution. Now 
suppose the state space of $\yv$ and $\av$ is
linearly transformed by a nonsingular matrix $\Sm$ as 
\begin{eqnarray}
     \yyv & = & \Sm\yv \\
     \aav & = & \Sm\av
\end{eqnarray}
where the transformation $\Sm$ could, for example, define 
tangent-space coordinates for a robot moving on a non-planar manifold
or, similarly, for gradient descent in a given function space.
Eq.(\ref{linmod}) can then 
be rewritten in the new coordinates as
\begin{equation}
     \yyv ~ = ~ \Sm\yv  ~ = ~ (\Sm\Amn\Smi)\cdot\Sm\av
     ~ = ~ (\Sm\Amn\Smi)\cdot\aav
\end{equation}
but for which  
\begin{eqnarray}
      \Sm\av  &  \neq & \pinv{(\Sm\Amn\Smi)}\cdot\yyv \\
      \Sm\av  &  \neq & \ginv{(\Sm\Amn\Smi)}\cdot\yyv 
\end{eqnarray}
In other words, the solution in the transformed space is not
the same as the original solution transformed to the new
space -- a situation that represents failure of a fundamental
sanity check. This is because the MP inverse and UC inverse
do not satisfy similarity consistency. 

Although the similarity-consistency (SC) condition of Eq.\!~(\ref{dinvcnst}) can be 
satisfied using the {\em Drazin}  inverse~\cite{drazin}, $\Adinv$, its use is not an option
for most engineering-related systems. Specifically, although the Drazin inverse
$\Adinv$ does have some practical applications~(\cite{cmr76,drazinApp}), 
its use is severely limited because $\Adinv$ may
have rank less than that of $\Amn$, which is often problematic
because it implies subspace information can be progressively
and irretrievably lost. In this paper we address this limitation of
the Drazin inverse by defining a new generalized inverse which satisfies
the condition of Eq.\!~(\ref{dinvcnst}) while also preserving
the rank of the original matrix.

\section{The SC Generalized Matrix Inverse}

One way to derive a generalized inverse that satisfies
a particular consistency condition is to identify a normal
form relative to that condition. For example, the singular
value decomposition (SVD) decomposes an arbitrary 
matrix $\Amn$ as 
\begin{equation}
      \Amn ~=~ \Um\Dm\Vmi
\end{equation}
where $\Um$ and $\Vm$ are unitary and $\Dm$ is a nonnegative
diagonal matrix of the singular values of $\Amn$. Critically, $\Dm$
is unique up to permutation, so all matrices with the same set of
singular values are equivalent up to left and right unitary
transformations. This allows the MP inverse to be derived as
\begin{equation}
   \Apinv ~ = ~ \Vm\tinv{\Dm}\Umi
\end{equation}
where $\tinv{\Dm}$ simply inverts the nonzero elements of $\Dm$.
From this the UC inverse can be similarly derived from a decomposition of the form
\begin{equation}
   \Amn ~ = ~ \Dm\Sm\Em
\end{equation}
where $\Dm$ and $\Em$ are nonnegative diagonal
matrices and $\Sm$ is unique up to left and right transformations
by unitary diagonal matrices~\cite{uhlmann}. The UC inverse can then be
obtained as
\begin{equation}
   \Aginv ~ = ~ \Emi\pinv{\Sm}\Dmi
\end{equation}
where the MP inverse can be applied because $\Sm$ is 
unique up to unitary diagonal transformations and
the MP inverse is consistent with respect to general
unitary transformations. 

Following a similar strategy, it can be recognized that the 
Jordan normal form of a square matrix $\Amn$ is
\begin{equation}
   \Amn ~ = ~ \Pm\Jm\inv{\Pm}
\end{equation}
where $\Pm$ is an arbitrary nonsingular matrix of the same
size as $\Amn$, and $\Jm$ is block-diagonal and unique up
to permutation of the blocks. Thus we can derive from this
decomposition the following similarity-consistent (SC) inverse:
\begin{equation}
   \Asinv ~ = ~ \Pm\pinv{\Jm}\inv{\Pm}
\end{equation}
where the MP inverse can be applied because $\Jm$ is 
unique up to block-diagonal permutation transformations,
which are of course unitary and thus permit the MP
inverse to be consistently applied.

What remains is to formally show that this definition of
$\Asinv$ is in fact consistent with respect to similarity
transformations, i.e., 
\begin{equation}
     \sinv{(\Sm\Amn\Smi)} ~ = ~ \Sm\Asinv\Smi .
\end{equation}
The Jordan decomposition of $\Amn$ is
\begin{equation}
   \Amn ~ = ~ \Pm\Jm\inv{\Pm}
\end{equation}
so the Jordan decomposition of $\Sm\Amn\Smi$ must
be
\begin{eqnarray}
   \Sm\Amn\Smi & = & \Sm(\Pm\Jm\inv{\Pm})\Smi\\
   ~ & = & (\Sm\Pm)\Jm(\inv{\Pm}\Smi)\\
   ~ & = & (\Sm\Pm)\Jm\inv{(\Sm\Pm)}
\end{eqnarray}
and thus the SC inverse of $\Sm\Amn\Smi$ is
\begin{eqnarray}
   \sinv{(\Sm\Amn\Smi)} & = & \sinv{\left((\Sm\Pm)\Jm\inv{(\Sm\Pm)}\right)}\\
   ~ & = &  (\Sm\Pm)\pinv{\Jm}(\inv{\Pm}\Smi)\\
   ~ & = & \Sm(\Pm\pinv{\Jm}\inv{\Pm})\Smi\\
   ~ & = & \Sm\Asinv\Smi.
\end{eqnarray}
This establishes similarity consistency, and because the Moore-Penrose
inverse preserves rank, the SC inverse does as well. The key difference
between the SC inverse and the Drazin inverse is that the latter 
enjoys the analytically important property that:
\begin{equation}
   \Amn\Adinv ~=~ \Adinv\Amn
\end{equation}
whereas product commutativity of $\Amn$ and $\Asinv$ 
is not a property of the SC inverse because $\pinv{\Jm}\Jm$
does not generally equal $\Jm\pinv{\Jm}$.

\section{Example and Practical Considerations}

Here we consider an illustrative numerical example
involving a nilpotent matrix
\begin{equation}
   \Amn ~=~ \left[\!\begin{tabular}{rrr}
                             4 & -1 & 2\\
                             7 & -2 & 3\\
                             -4 & 1 & -2
                     \end{tabular}\right]
\end{equation}
with Jordan form $\Amn=\Pm\Jm\inv{\Pm}$ where
\begin{equation}
   \Pm ~=~ \left[\!\begin{tabular}{rrr}
                             1 & -2 & -3\\
                             2 & -5 & 10\\
                            -1 & 2 & -2
                     \end{tabular}\right]
~ \mbox{and} ~ ~
   \Jm ~=~ \left[\!\begin{tabular}{rrr}
                             0 & 1 & 0\\
                             0 & 0 & 1\\
                             0 & 0 &  0
                     \end{tabular}\right].
\end{equation}
It can be verified that 
\begin{equation}
   \Asinv ~=~ \left[\!\begin{tabular}{rrr}
                             -2 & 1 & 2\\
                             10 & 0 & 15\\
                             8 & -2 & 2
                     \end{tabular}\right]
~ \mbox{and} ~ ~
   \Adinv ~=~ \left[\!\begin{tabular}{rrr}
                             0 & 0 & 0\\
                             0 & 0 & 0\\
                             0 & 0 &  0
                     \end{tabular}\right].
\end{equation}
where the Drazin inverse of a nilpotent matrix is always
zero, which graphically demonstrates complete loss of 
rank\footnote{The integer matrices in this section were constructed
using the method of~\cite{hanson} so that inverses are also 
integer matrices. This allows the results to be more conveniently
verified by direct hand calculations.}.
We now consider a matrix $\Mm=\Sm\Amn\Smi$ similar to $\Amn$
with
\begin{equation}
   \Sm ~=~ \left[\!\begin{tabular}{rrr}
                             -1 & 4 & 0\\
                             -2 & 9 & 3\\
                              1 & -4 & 1
                     \end{tabular}\right].
\end{equation}
Similarity consistency for the Drazin inverse is trivially
maintained because 
$\dinv{\Mm}=\dinv{(\Sm\Amn\inv{\Sm})}=\Sm\dinv{\Amn}\inv{\Sm}={\bf 0}$.
By contrast, it can be verified that the rank of the SC inverse
\begin{equation}
\sinv{\Mm} ~=~ \left[\!\begin{tabular}{rrr}
                             -819 & 167 & -443\\
                             -2301 & 464 & -1255\\
                              663 & -137 & 355
                     \end{tabular}\right]
\end{equation}
is the same as $\Mm$ and that the similarity consistency condition
$\sinv{(\Sm\Amn\inv{\Sm})}=\Sm\sinv{\Amn}\inv{\Sm}=\sinv{\Mm}$
is satisfied. And of couse it can also be verified
that $\sinv{\Mm}$ satisfies the generalized inverse
conditions $\Mm\sinv{\Mm}\Mm=\Mm$ and $\sinv{\Mm}\Mm\sinv{\Mm}=\sinv{\Mm}$.
When this example is numerically evaluated using the Matlab {\sf jordan(M)} 
and {\sf pinv(M)} functions the SC inverse is found to be accurate to machine precision.
However, when applied in two different high data-rate robotic control applications 
the observed errors using the SC inverse were sufficiently large to cause the 
control process to diverge rapidly\footnote{These simulation experiments 
were performed using variations on models for a robot arm and a 
rover examined in~\cite{zhang}. Unit consistency was the focus of that
work and the two failed tests of the SC inverse were neither analyzed nor reported.}. 

The numerical instability of the Jordan decomposition is well known~\cite{golub}, 
and this problem is exacerbated when applied with the MP inverse because large errors 
can be introduced into singular values that should be zero, which can lead to huge
spike errors when inverted. It is not clear to what 
extent this issue can be mitigated when these two operations are used in
succession to compute the SC inverse. In many practical dynamical-system applications 
the system matrix changes relatively slowly over time, so it may be possible
to track the evolution of its distribution of singular values to actively avoid 
catastrophic discontinuities in the case of a singular value near the 
threshold between zero and nonzero when the MP inverse is applied
during evaluation of the SC inverse. More specifically, a tracked presumed-zero 
singular value that drifts above the MP threshold can be maintained as zero while 
a tracked presumed-nonzero singular value can be maintained as nonzero even if it 
drifts below the threshold.

A potentially more general approach for mitigating some of the numerical fragility is 
to either jointly solve the Jordan+MP calculation or to relax the form of the 
decomposition, e.g., so that the central matrix is unique only up to orthonormal
similarity to a complex symmetric matrix rather than
permutation similarity to a Jordan matrix. For example, any matrix $\Amn$
can be decomposed\!~\cite{HJ1} as
\begin{equation}
   \Amn ~ = ~ \Sm\Xm\Smi  \label{symdecomp}
\end{equation}
where $\Xm$ is complex symmetric and thus
the SC inverse can be defined as:
\begin{equation}
   \Asinv ~ = ~ \Sm\pinv{\Xm}\Smi   \label{symsc}.
\end{equation}
That the decomposition of
Eq.\!~(\ref{symsc}) is only unique up to unitary/orthonormal 
similarity can be inferred from the fact that for
complex symmetric $\Ym$ and orthonormal $\Qm$ ($\inv{\Qm}=\Qmt$),
$\Ym=\Qm\Xm\Qmt$ implies
\begin{eqnarray}
   \Asinv & = & \Sm\pinv{\Xm}\Smi\\
    ~ & = & \Sm\pinv{(\Qm\Ym\Qmt)}\Smi \\
     ~ & = &      (\Sm\Qm)\pinv{\Ym}\inv{(\Sm\Qm)} 
\end{eqnarray}
where the last step is obtained by exploiting the unitary consistency property of the MP inverse.
It is reasonable to expect that this less rigid definition will be more amenable
to numerically-robust evaluation than is possible using a construction requiring
an explicit Jordan decomposition.

\section{Fixed-Rank Assumption}

In some applications it is known that the rank of the
time-evolving system matrix can potentially change, 
e.g., during transitions through gimbal-lock 
configurations, while in others
it can be expected to remain constant. In the latter case
a fixed-rank assumption can be exploited to mitigate
the problem of spurious singular values when applying
the MP inverse to the Jordan matrix during evaluation of
the SC inverse. Specifically, if a fixed rank of $k$ is 
known/assumed then the MP inverse can be evaluated
explicitly using the SVD with all but the $k$
largest singular values set identically to zero instead of having
to discriminate with a threshold.  

The efficacy of exploiting the fixed-rank assumption 
can be demonstrated in a carefully controlled experiment
involving a dynamical 3-dimensional system that is
confined to a 2-dimensional subspace such that the
$3\times 3$ system matrix can be assumed to have a fixed
rank of $2$. The testing model can be further tailored 
so that the system matrix is guaranteed to never be
diagonalizable and thus is maximally susceptible to the 
numerical issues discussed in the previous section. We
achieve this by defining the Jordan form of the system
as a function of time as
\begin{equation}
   (\Qm(t)\Sm)\,\Jm\,\inv{(\Qm(t)\Sm)}
\end{equation}
with
\begin{equation}
\Jm = \begin{array}{|ccc|}
             1 & 1 & 0\\
             0 & 1 & 0\\
             0 & 0 & 0
          \end{array}
~,~ ~
\Qm(t) = \begin{array}{|ccc|}
                    \cos(t) & -\sin(t) & 0\\
                    \sin(t)  & ~\cos(t) & 0\\ 
                    0 & 0 & 1
               \end{array}
\end{equation}
and $\Sm$ is an arbitrary fixed nonsingular matrix. Thus
the system matrix is singular and has an orbit period of $2\pi$ 
during which it is never diagonalizable. The plot of Figure\!~1 
shows the maximum absolute elementwise deviation between the 
computed SC inverse and the ground-truth known SC inverse 
for each of 300 uniform timesteps $t$, $0\leq t<2\pi$. It
is seen that spike errors occur at a significant fraction of the 
timesteps. These spikes are actually of magnitude $>\!10^{14}$ 
and result from inverting spurious singular values 
that exceed {\sf pinv}'s zero threshold. In tests of many
randomly-generated $\Sm$ matrices the number of
spike errors over the entire 300-timestep simulation ranged from
zero to more than $50\%$, with Figure\!~1 showing
an example of the latter. 
\begin{figure}
	\centering
	\includegraphics[width=3.5in]{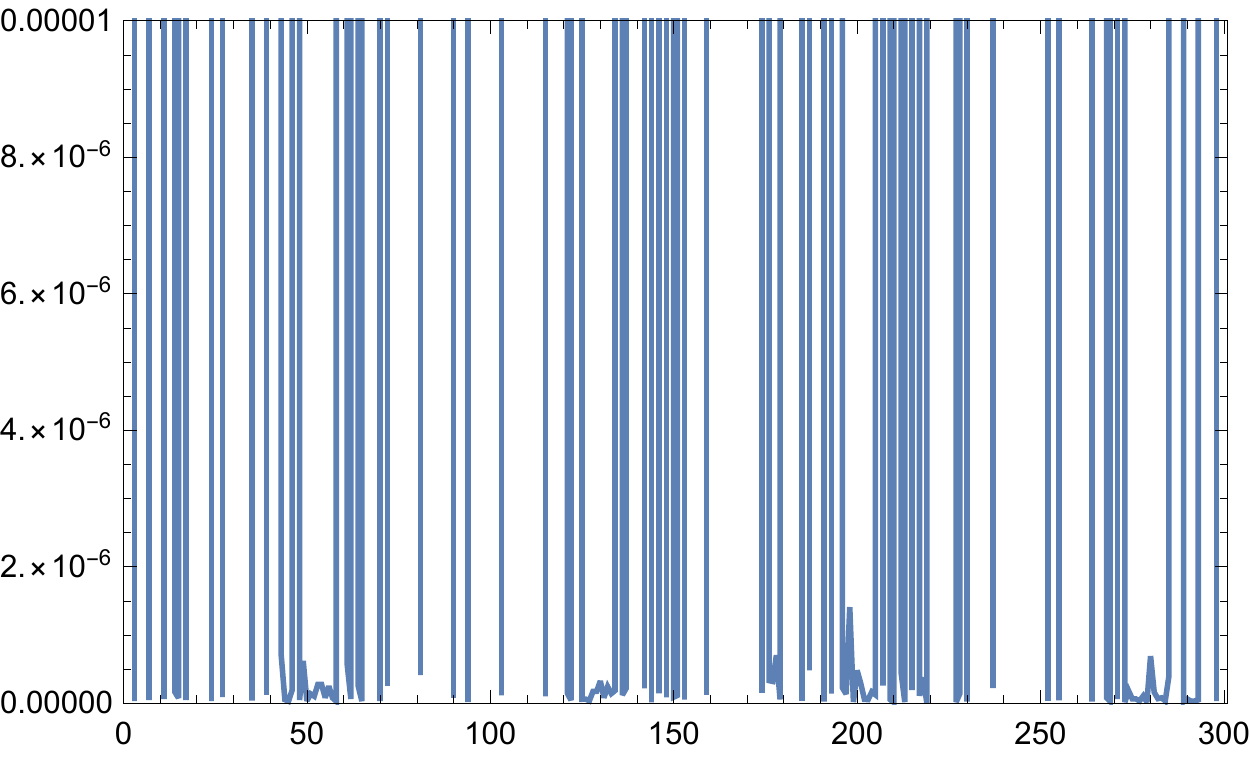}
	\caption{Error magnitudes (vertical axis) resulting from standard evaluation of the
                         SC inverse from the Jordan form over 300 timesteps (horizontal axis).
                         Clipped magnitudes are actually of size $\approx 10^{14}$.}        
\end{figure}

Figure\!~2 shows the results of the same
scenario but with the smallest singular value of the Jordan matrix
always treated as identically zero for computation of the SC inverse at
each timestep. The results shown in Figure\!~2 are typical of those
observed for the fixed-rank method in the battery of tests over many
randomly-generated $\Sm$ matrices. By design it is likely that the model 
used for these tests overestimates the frequency of spike errors that can be 
expected from standard evaluation of the SC inverse in real-world 
applications, but because any such errors will be time-correlated it likely will
not be practical to simply identify and discard them (e.g., apply some sort of 
heuristic outlier filter) on the assumption that they will only occur sporadically.
In other words, the intrinsic numerical instability of the SC inverse cannot be
mitigated in the general case. However, the results of Figure\!~2 suggest
that the problem may be effectively mitigated in a large class of applications.
\begin{figure}
	\centering
	\includegraphics[width=3.5in]{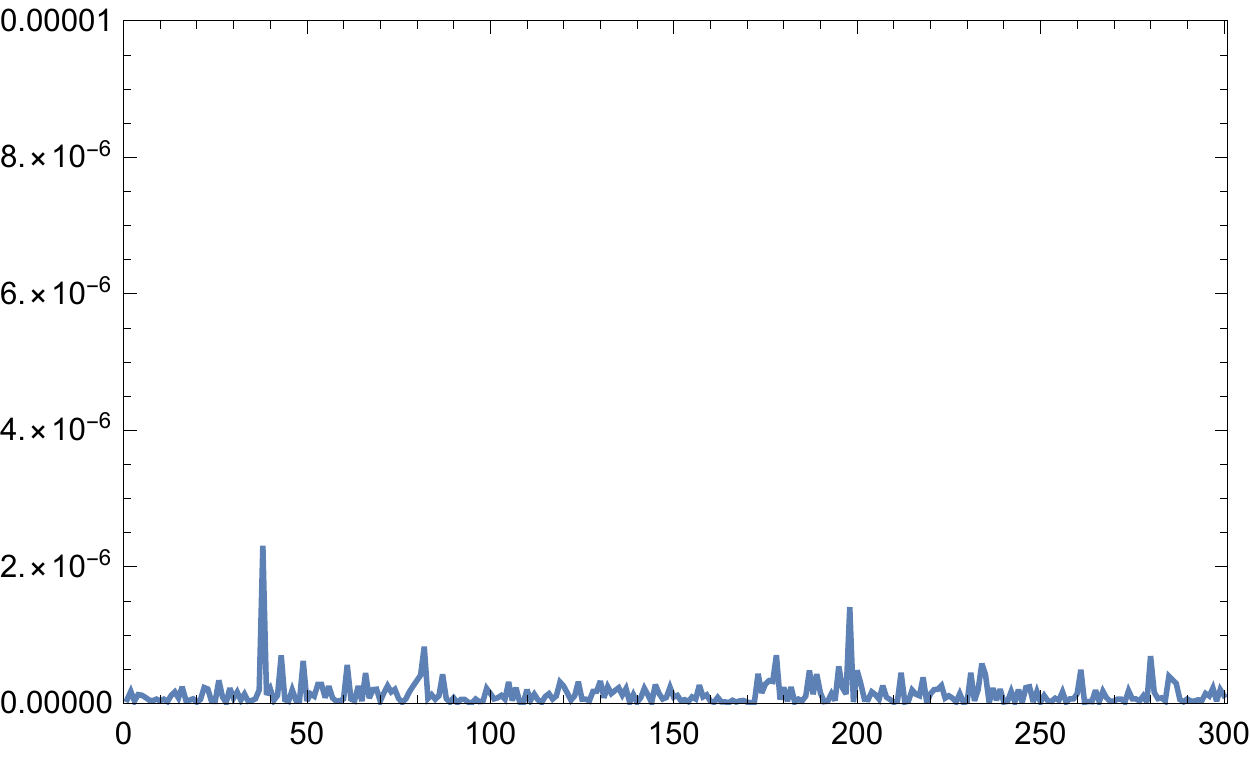}
	\caption{Error magnitudes (vertical axis) resulting from fixed-rank evaluation of the
                         SC inverse from the Jordan form over 300 timesteps (horizontal axis). The mean
                         error magnitude is $1.4420e$-$07$ and the maximum magnitude is $2.3101e$-$06$.}        
\end{figure}

\section{Discussion}

Similarity transformations capture the most general linear transformations that can
be applied to the common coordinate frame of a set of input and output 
variables, e.g., as can arise in machine learning applications in which little is 
known about the best space in which to represent an input-output mapping.
In this paper we defined a new similarity-consistent generalized matrix
inverse that preserves/maintains the rank of the original matrix, a property
which is not preserved by the Drazin inverse. It has been noted that the
numerical evaluation of the SC inverse can be a challenging practical
limitation. One possible approach for addressing this problem is to 
relax the structure of the decomposition from requiring similarity to a Jordan
matrix to similarity to a complex symmetric matrix for which uniqueness
is only required up to orthonormal similarity, rather than permutation 
similarity, so the MP inverse is still applicable. This less-rigid
formulation may be more amenable to numerically-robust evaluation,
especially in conjunction with the fixed-rank method.
Regardless of the choice of decomposition, 
it would be interesting to examine whether recent
methods developed for evaluating the Drazin inverse can be applied
to the SC inverse~\cite{drazinalg,pan}. 

Despite its computational challenges, the SC inverse is the only 
general option available for applications that demand consistency with
respect to arbitrary linear transformations of a common coordinate
frame. We have provided evidence that these challenges can be 
effectively addressed in applications in which the rank of the system
can be assumed fixed during its time evolution, but we have also
emphasized that intrinsic numerical instabilities represent a 
fundamental obstacle in the general case. 
This motivates continued work toward identifying special properties
of particular applications that can be exploited to permit the SC
inverse to be practically applied.

\end{document}